\documentclass[12pt]{amsart}
\usepackage{amsmath, amsthm, amsfonts}

\usepackage{caption}
\usepackage{subfig}

\usepackage{graphicx}
\usepackage{hyperref}
\usepackage{float}

\newtheorem{theorem}{Theorem}

\DeclareMathOperator{\diam}{diam}

\theoremstyle{definition}
\newtheorem{definition}{Definition}

\usepackage{mathtools}

\newcommand{\eqdef}{=\vcentcolon}

\newcommand{\R}{\ensuremath{\mathbb{R}}}

\newcommand{\Z}{\ensuremath{\mathbb{Z}}}

\newcommand{\N}{\ensuremath{\mathbb{N}}}

\begin{document}

\title[Brownian Earthworm]{Hausdorff Dimension of the Brownian Earthworm\\
(WXML Winter 2017 Report)}
\author[]{Clayton Barnes,  Krzysztof Burdzy,  Alex Forney, \\
Zach Dingels, Xiyi Yan, and Ran Zhao}
\address{Department of Mathematics, University of Washington}
%\date{Winter 2017}
\maketitle

\section{Introduction}

\subsection{Notions of Dimension}
When we refer to the ``dimension'' of a set, we are referring to some notion of the set's size. In particular, we are addressing two common notions of dimension for a set $A \subset \R^n$:
\begin{itemize}
\item Minkowski dimension, and
\item Hausdorff dimension.
\end{itemize}
The Minkowski dimension (sometimes called Minkwski-Bouligand dimension) is often referred to simply as box-counting dimension. As a consequence of the finite capacity of computing power, one often uses the box-counting dimension as an approximation of the Hausdorff dimension. A more thorough treatise on box-counting and Hausdorff dimension, as well as other notions of dimension, can be found in \cite{falconer}.

\subsubsection{Box-Counting Dimension}
Given a set $A$ and $n \in \N$, let $N(\varepsilon)$ be the minimum number of squares (or balls) of side length (resp. radius) $\varepsilon$ which can cover $A$.

Consider the rate at which $N(\varepsilon)$ grows as $\varepsilon \to 0$. That is, take
\begin{equation} 
d = \lim_{\varepsilon \to 0} \frac{\log{N(\varepsilon)}}{\log{1/\varepsilon}}. 
\label{eq: BoxCounting}
\end{equation}
Provided this limit exists, we call $d \eqdef \dim_{B}(A)$ the box-counting dimension of $A$.

For example, let $A = [0, 1] \subset \R$. Given a square with sides of length $\varepsilon_n = 1/n$, $n$ squares are needed to cover $A$. Hence, $N(\varepsilon_n) = n$. By computing the limit given in \eqref{eq: BoxCounting}, we find that $d = 1$. This implies that the box-counting dimension of the line segment is $1$.

Similarly, consider the case where $A = [0, 1] \times [0, 1] \subset \R^2$. Now, it would take $n^2$ squares with sides of length $1/n$ to cover $A$, which implies that $N(n) = n^2$. The limit in \eqref{eq: BoxCounting} gives that $d = 2$.

These two examples correspond with the intuitive understanding that a line is 
one-dimensional and that a square is two-dimensional.

\subsubsection{Hausdorff Dimension}
In order to define the Hausdorff dimension, we must first define the Hausdorff measure.

\begin{definition}
Let $U \subseteq \R^n$. We say that $\diam{U}$ is the \textbf{diameter} of $U$, where
\[ \diam{U} = \sup\{\|x - y\| \mid x, y \in U\}. \]
\label{def:diam}
\end{definition}

\begin{definition}
Let $s \geq 0$, $A \subset \R^n$, and $\delta > 0$. Define
\[ H_{\delta}^{s} (A) = \inf{\left\{ \sum_{k=1}^{\infty} (\diam{U_k})^s \mid \bigcup_{k=1}^{\infty} U_k \supseteq A \text{ and } \diam{U_k} < \delta \right\}}. \]
Further, define
\[ H^{s} (A) = \lim_{\delta \to 0} H_{\delta}^{s} (A). \]
We call $H^{s} (A)$ the \textbf{$s$-dimensional Hausdorff measure} of $A$.
\end{definition}

\noindent Now, let $A \subset \R^n$, $0 < \delta < 1$, and $0 \leq s < t$. Suppose we have
\[ \bigcup_{k=1}^{\infty} U_k \supseteq A. \]
It follows that
\[ \sum_{k=1}^{\infty} (\diam{U_k})^t \leq \delta^{t-s} \sum_{k=1}^{\infty} (\diam{U_k})^s, \]
which implies that
\[ H_{\delta}^{t} (A) \leq \delta^{t-s} H_{\delta}^{s} (A). \]
Taking the limit of this inequality as $\delta \to 0$, we can see that if $H_{\delta}^{s}$ is not infinity, then $H_{\delta}^{t} (A) = 0$. In particular, the point at which $H_{\delta}^{s}$ goes from infinity to $0$ is the Hausdorff dimension of the set $A$. This phenomenon can be seen in the following plot:

\begin{figure}[H]
\centering
\includegraphics[scale = .4]{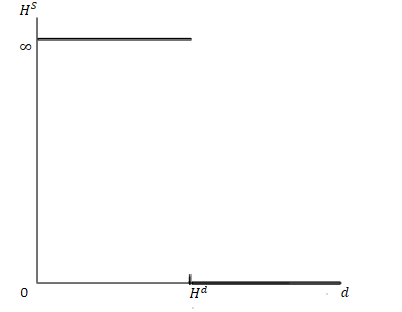}
\end{figure}

\noindent More precisely, we have
\begin{definition}
The \textbf{Hausdorff dimension} of $A \subset R^n$ is given by
\[ \dim_{H} (A) = \inf{ \left\{ s \geq 0 \mid H^{s} (A) = 0 \right\}} = \sup{ \left\{ s \geq 0 \mid H^{s} (A) = \infty \right\}}. \]
\end{definition}

\noindent We have a comparison with the box-counting dimension by the following
\begin{theorem}(Proposition 3.8, \cite{falconer})
For any set $A,$
$$\dim_{H}(A) \leq \dim_{B}(A).$$
\end{theorem}

\subsection{Brownian Frontier}

\ \ 

\noindent\begin{minipage}{0.3\textwidth}
\includegraphics[width=\linewidth]{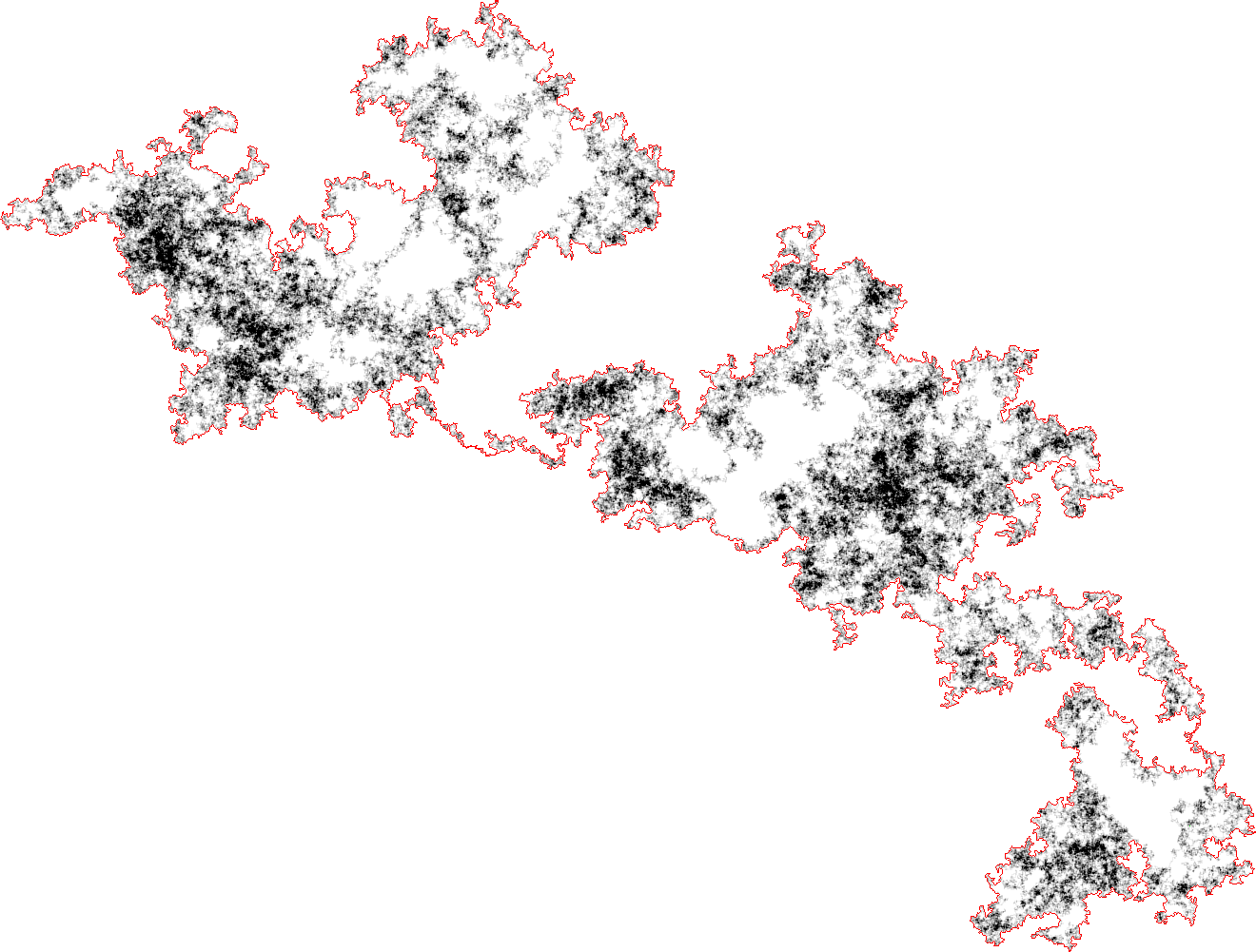}
\end{minipage}
\hfill
\begin{minipage}{0.3\textwidth}
\includegraphics[width=\linewidth]{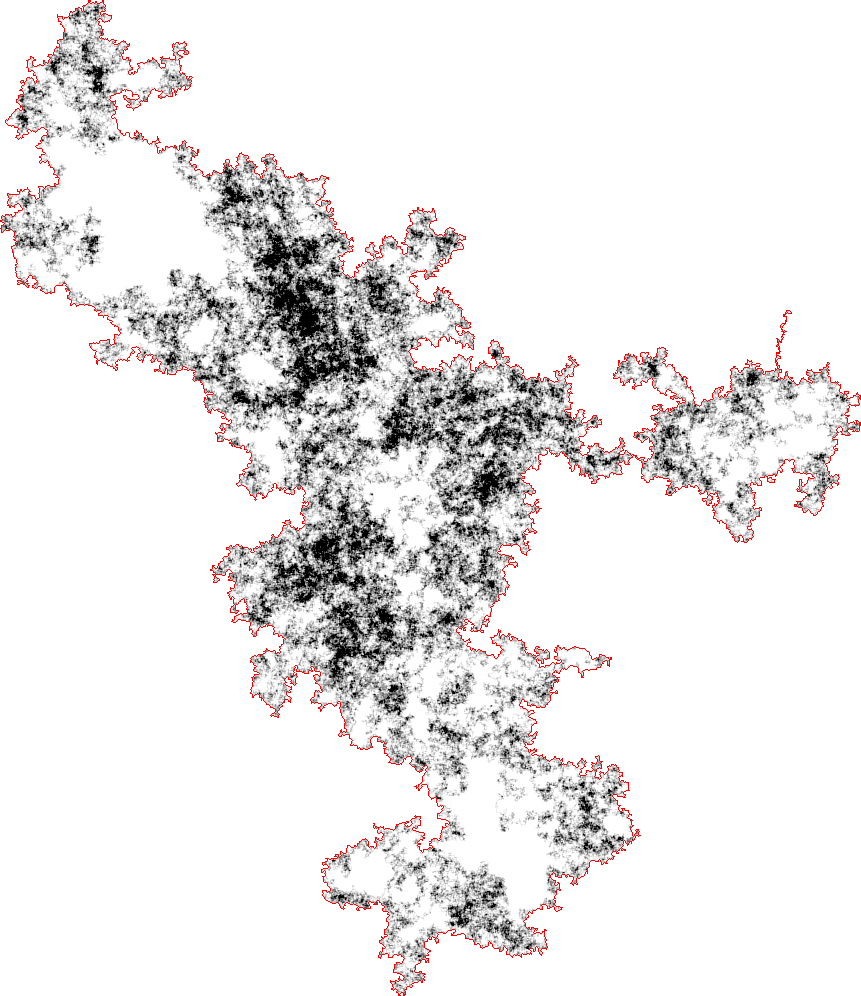}
\end{minipage}
\hfill
\begin{minipage}{0.3\textwidth}
\includegraphics[width=\linewidth]{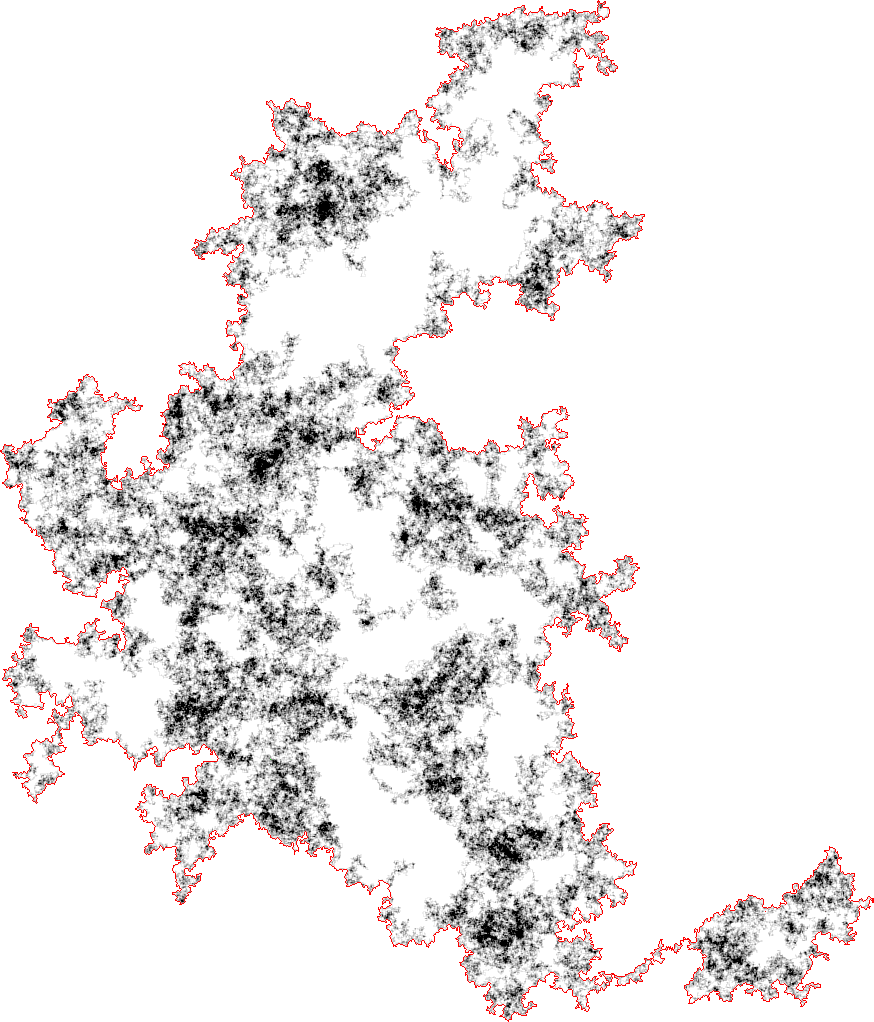}
\end{minipage}
\vspace{1em}

\noindent The first objective of this project was to analyze Brownian motion. Specifically, we wanted to calculate the Hausdorff dimension of the Brownian frontier. For our purposes of discrete approximation, we appeal to Donsker's theorem, which says that Brownian motion is a limit of scaled random walks. The existence of Brownian motion and Donsker's invariance principle are beyond our scope, but a reference can be found in \cite{morters2010brownian}.

\begin{definition}
A one-dimensional Brownian motion $\{B(t) \mid 0 \leq t < \infty\}$ is a stochastic
process with the following properties:
\begin{itemize}
\item $B(t)$ is continuous almost surely.
\item For any $0 \leq s < t$, $B(t) - B(s) \overset{d}{=} N(0, t-s)$.
\item For any $0 \leq t_1 < t_2 < \cdots < t_k$, $X_{t_1}, X_{t_2} - X_{t_1}, \dots,
X_{t_k} - X_{t_{k-1}}$ are all independent.
\end{itemize}
\end{definition}

\noindent A two-dimensional Brownian motion is simply a pair $(B_1, B_2)$ where $B_1$ and $B_2$ are two independent one-dimensional Brownian motions.

Let $X_i, i \in \N$ be independent and identically distributed random variables $P(X_i = 1) = 1/2 = P(X_i = -1)$. Define $\displaystyle S_n = \sum_{i = 1}^n X_i$, and $S_n(t) = S_{\lceil nt \rceil}/\sqrt{n}$. Then,

\begin{theorem}(Donsker's Theorem)
$\{S(t) \mid t \in [0, 1]\} \overset{d}{\to} \{B(t) \mid t \in [0, 1]\}$ as distributions on the space of continuous real valued functions on $[0, 1]$.
\end{theorem}

Further, we have

\begin{definition} The \textbf{frontier} of Brownian motion is the boundary of the unbounded component of the complement of Brownian motion.
\end{definition}

\noindent In the above images, the frontier is highlighted in red.

\subsection{Brownian Earthworm}

\ \ 

\noindent\begin{minipage}{0.3\textwidth}
\includegraphics[width=\linewidth]{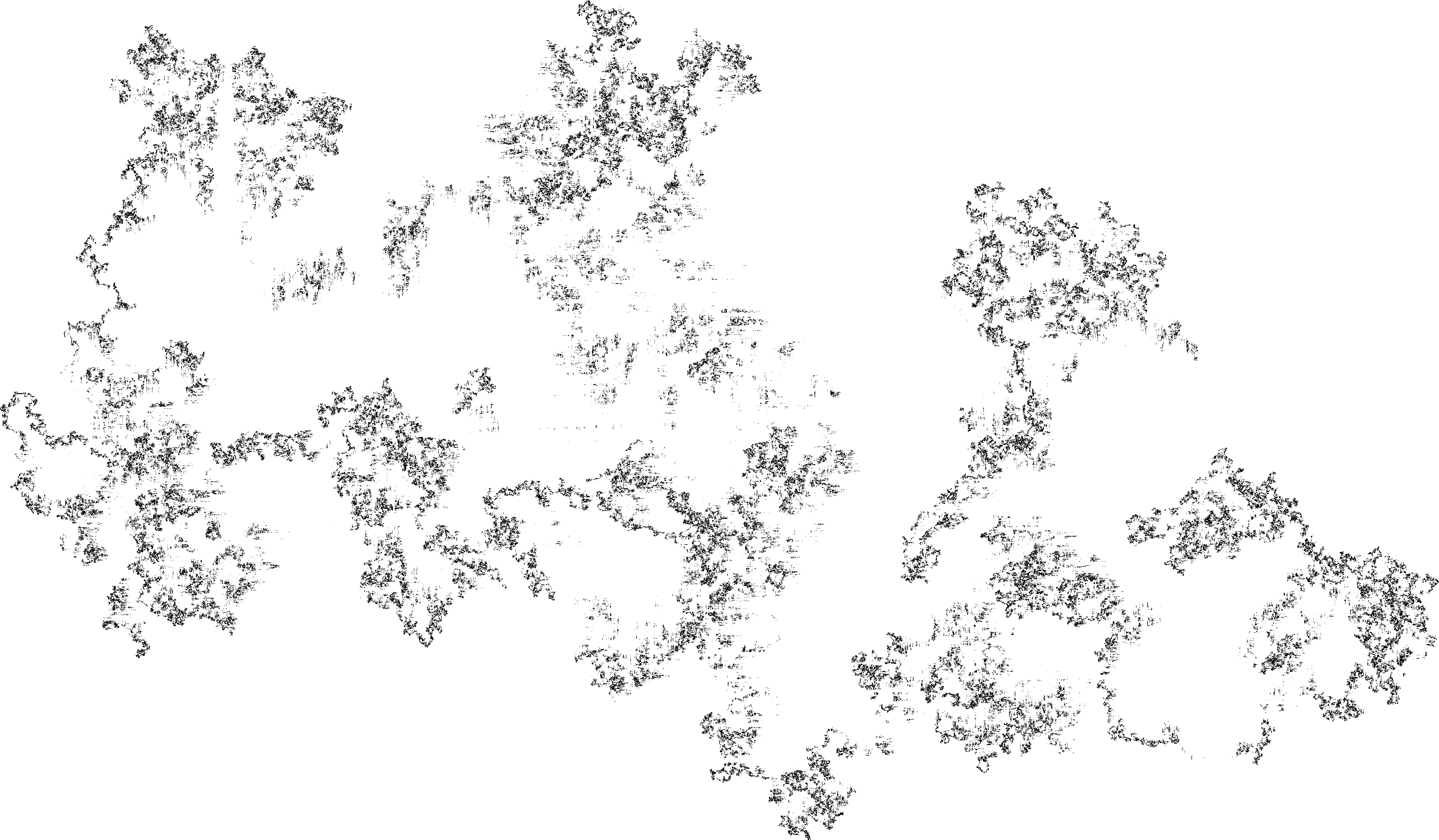}
\end{minipage}
\hfill
\begin{minipage}{0.3\textwidth}
\includegraphics[width=\linewidth]{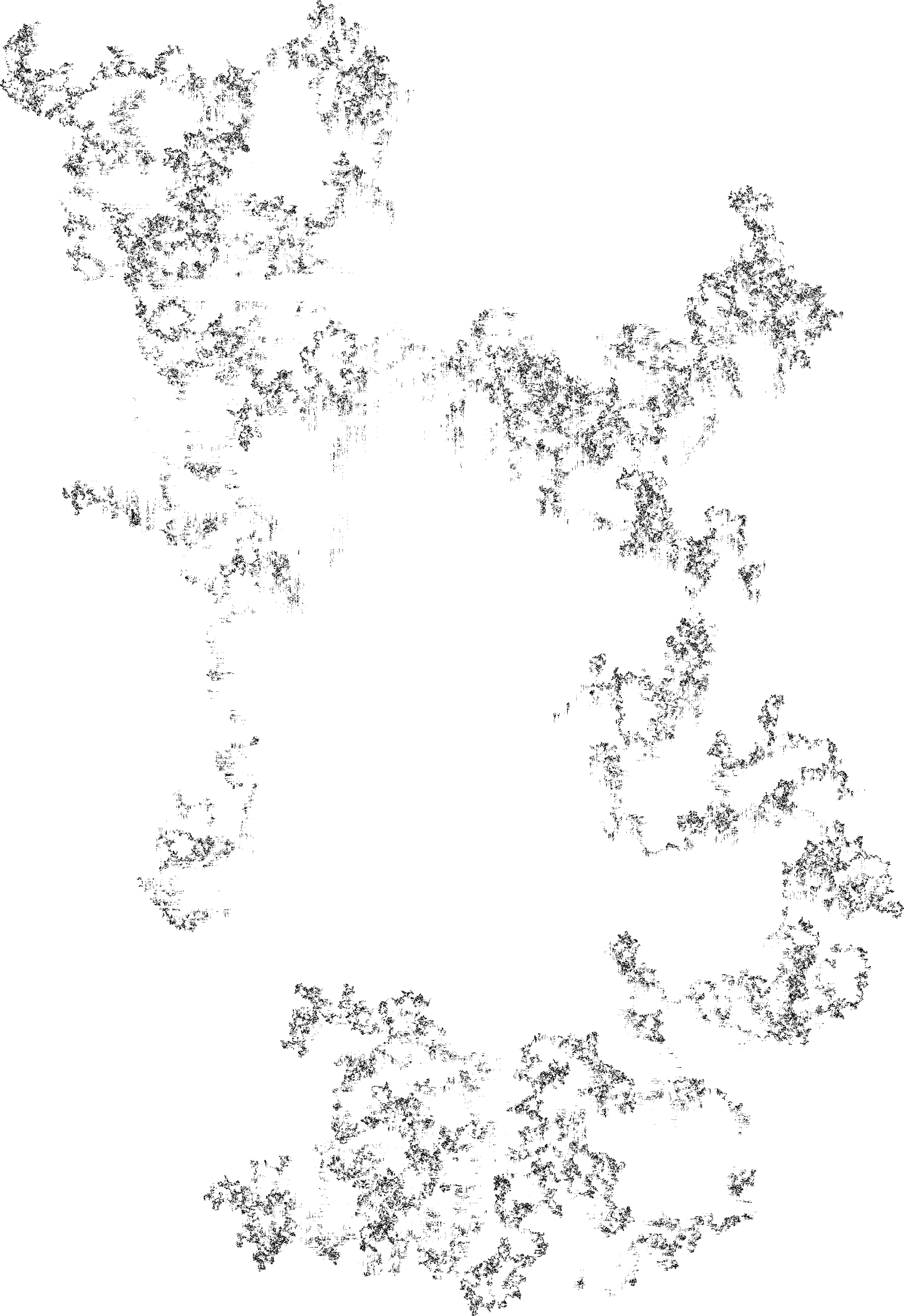}
\end{minipage}
\hfill
\begin{minipage}{0.3\textwidth}
\includegraphics[width=\linewidth]{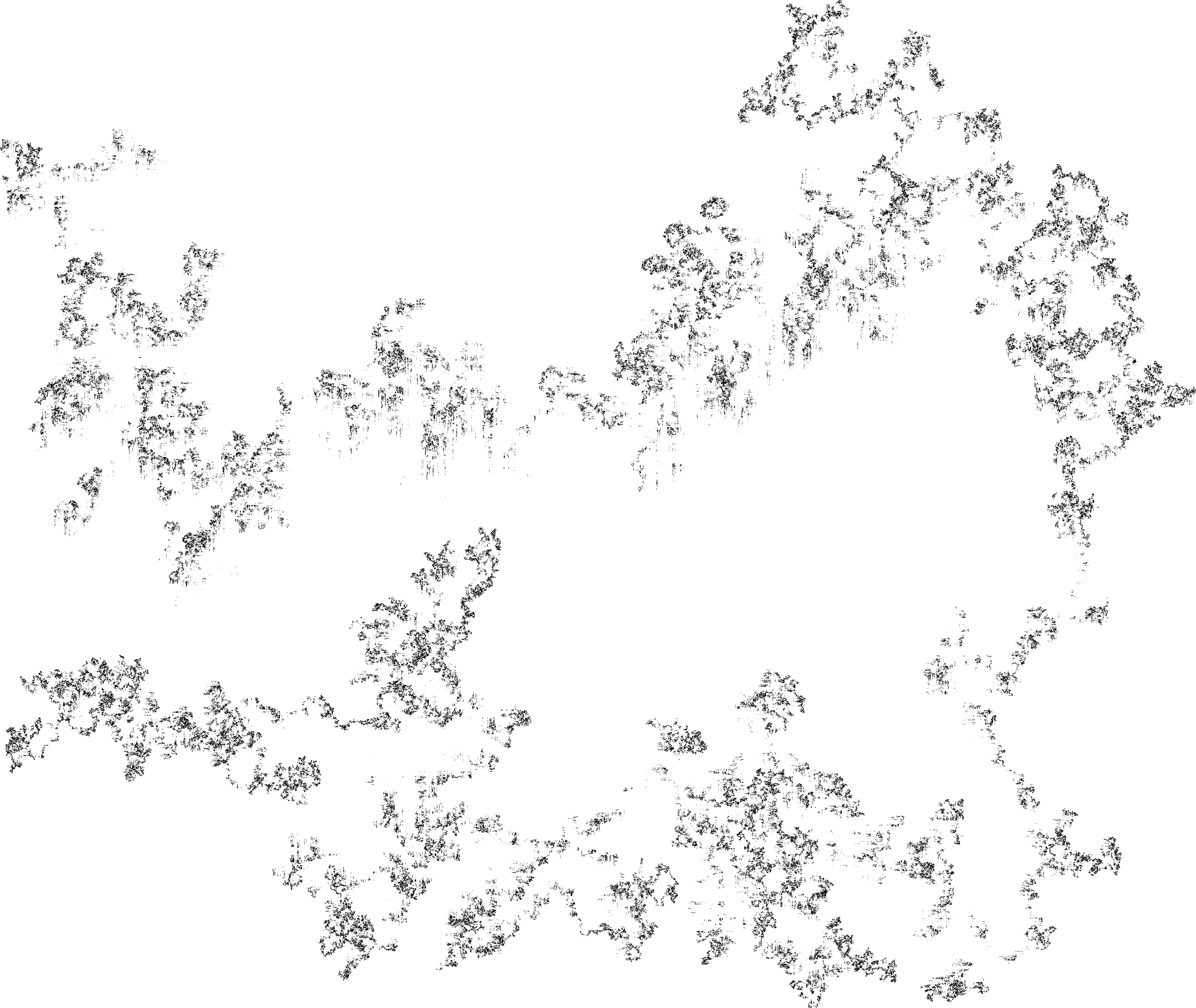}
\end{minipage}
\vspace{1em}

\noindent Suppose that there is a particle of soil at every point in $\Z^2$. The Brownian earthworm model \cite{earthworm} is a random walk in the plane with the addition that the earthworm leaves behind a set of holes. That is, the earthworm pushes the particles of soil in front of it and leaves a hole in its wake. If the earthworm moves in the direction of an already-created hole, then the displaced soil will fill in the hole. This model can be visualized by the above images, where each black point is a hole left behind by the earthworm.

\section{Methods}

\subsection{Computation}
Our first task was to simulate Brownian motion and isolate the frontier. So, we set up a simulation of a random walk in $\Z^2$ in Python and allowed it to run for millions of steps. Obviously, our simulation only took on discrete values, which meant that we had to scale our walks by some value so that they converged to Brownian motion (cf. Donsker's Theorem). In order to actually calculate the frontier points, which we found to be a nontrivial task, we contained our walk of fixed length in a rectangle of fixed size. Then, we took points from the boundary of the rectangle and moved them ``inward'' until they landed on a point in the walk. We repeated this process until the entire frontier was covered.

The simulation of the earthworm model again began with a random walk in $\Z^2$ in Python. Instead of letting the simulation run for a fixed number of steps and then performing calculations on our set of data points, we maintained a list of the set of holes as the simulation was running. This was largely because holes can be filled in as the earthworm travels and because we found a list to be more efficient than, say, a two-dimensional array representation of the lattice.

\subsection{Counting Method}
This method of estimating the Hausdorff dimension is directly motivated by the box-counting dimension approximation. It can be applied to any set in $\Z^2$, but we describe how it applies specifically to our models.

For a simulated random walk with $n$ steps, let $S \subset \Z^2$ be the set in which we are to estimate the Hausdorff dimension. For instance, $S$ could be the set of points in the frontier or the collection of holes from the Brownian earthworm. It follows that $S$ will have a diameter $\diam{S}$. As an aside, this diameter in either model is close to the maximum distance a random walker travels from the origin in $n$ steps, which has an expected value proportional to $\sqrt{n}$ which we see in the plots of the diameter.

Say we observe that $\diam{S}$ is approximately $n^d$ for some $d \in \R$ and similarly that the size of $S$ is approximately $n^h$ for some $h \in \R$. By the definition of diameter, $S$ is contained in a square of side length $n^d$. Further, there are $n^{2d}$ vertices in this square. As a result, each vertex in this square can be thought of as small box of side length $\varepsilon_n = 1/n^d$. (Essentially, this amounts to normalizing the large square to have unit mass.) This means that the number of boxes needed to cover $S$, call it $N(\varepsilon_n)$, is $|S| \approx n^h.$ By plugging these values directly into the corresponding entries for the box counting dimension given in \eqref{eq: BoxCounting}, we see that
\[ \dim_{B}{S} \approx \frac{\log N(\varepsilon_n)}{\log(1/\varepsilon_n)} = \frac{\log n^h}{\log n^d}. \]
Since $\varepsilon_n \to 0$ as $n\to \infty$, we estimate the box counting dimension by taking many observations for increasingly large values of $n$.

In practice, this amounts to plotting the growth rate of the size of the set against the growth rate of the diameter on a log-log scale.

\subsection{Averaging Method}
The method stems from the intuition behind the Hausdorff dimension. If a set $S$ has Hausdorff dimension $d$ and $s$ is a point in $S$, then, in layman's terms, the volume (Hausdorff measure) of a ball of radius $r$ around $s$ is approximately $r^d$. We can isolate $d$ by taking logarithms and dividing by $\log{r}$.

Since our set is discrete, we allow the number of steps to increase to infinity and let $r$ vary within a finite range. We then record the ``volume'' (number of vertices) within a radius of $r$ and further take the average volume across all values of $r$.

More precisely, let $r \in \N$. For each element $s$ in our set of unknown dimension, we want to find the number of elements that are within a ball of radius $r$ from $s$. That is, we want to find $|B_{s, r}| = |\{x \in S \mid |s - x| \leq r\}|$. Then, we calculate the average, $Q_{r}$, of all such $|B_{s, r}|$.
We see that the number $Q_{r} = r^\alpha$ for some $\alpha \in \R$. Again, we can find the dimension of the set by plotting $Q_{r}$ against $r$ on a log-log scale.

\section{Results}

\subsection{Brownian Frontier}
Since the dimension of the Brownian motion is known to be $2$ and the dimension of the Brownian frontier is known to be $4/3$ \cite{frontier}, we wished to verify that our two methods for calculating the Hausdorff dimension produced the desired results.

With the Counting Method, we first calculated the coefficient corresponding to the growth rate of the diameter to be $0.4986$ (Figure \ref{fig: diam} below). We then used this value to find the dimensions of Brownian motion (Figure \ref{fig: area}) and the Brownian frontier (Figure \ref{fig: frontier}) to be $1.8662$ and $1.3364$, respectively. Our result for the dimension of the frontier was close to the actual value, although our calculated dimension of Brownian motion was not as close as we would have liked.

\begin{figure}[H]
\begin{minipage}{\linewidth}
  \centering
  \includegraphics[width=0.4\linewidth]{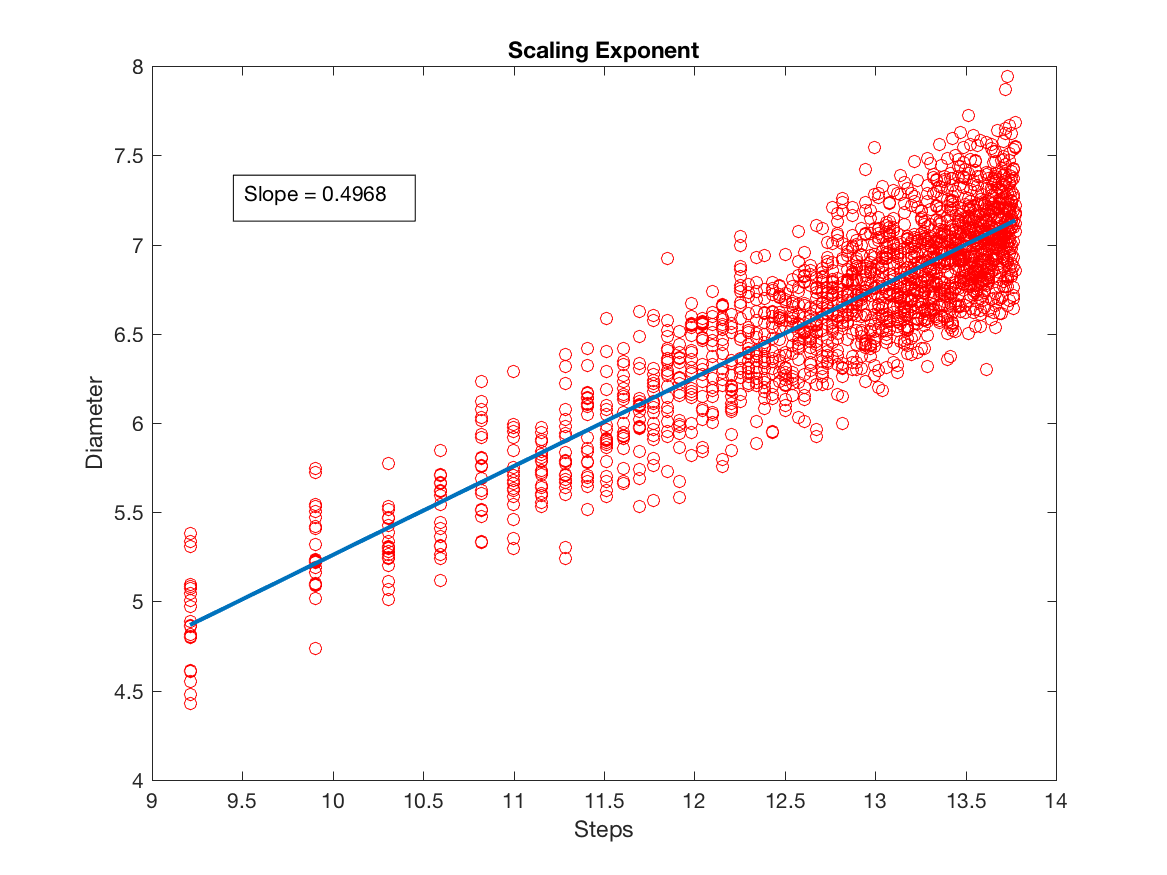}
  \caption{Diameter}
  \label{fig: diam}
\end{minipage}
\vspace{1em}
\begin{minipage}[b]{0.4\linewidth}
  \centering
  \includegraphics[width=\linewidth]{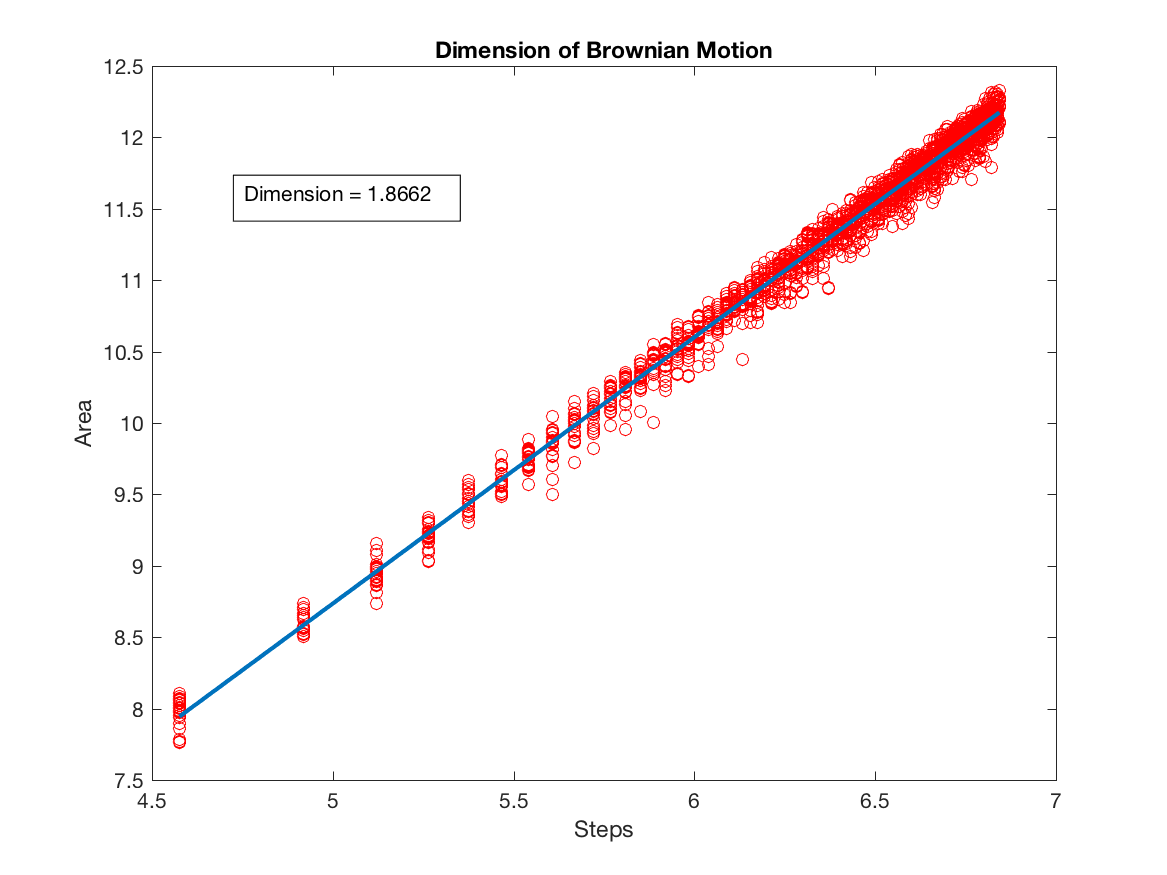}
  \caption{Brownian Motion}
  \label{fig: area}
\end{minipage}
\hfill
\begin{minipage}[b]{0.4\linewidth}
  \centering
  \includegraphics[width=\linewidth]{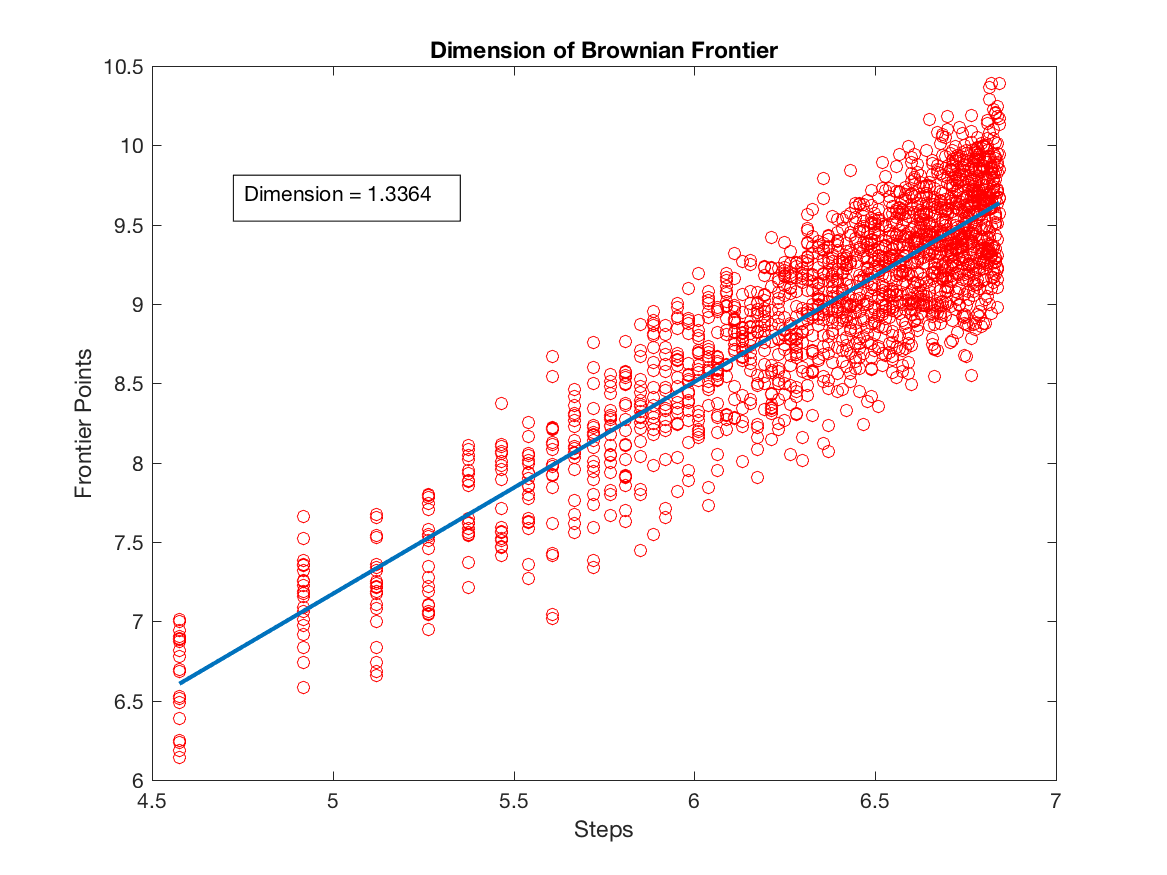}
  \caption{Brownian Frontier}
  \label{fig: frontier}
\end{minipage}
\end{figure}

We next used the Averaging Method to attempt to corroborate our results. We focused this method on the Brownian frontier in particular since this was our main interest. In fact, this second method gave a very similar value, $1.3430$, to our first calculated value and the actual value.

\begin{center}
\includegraphics[width=0.6\linewidth]{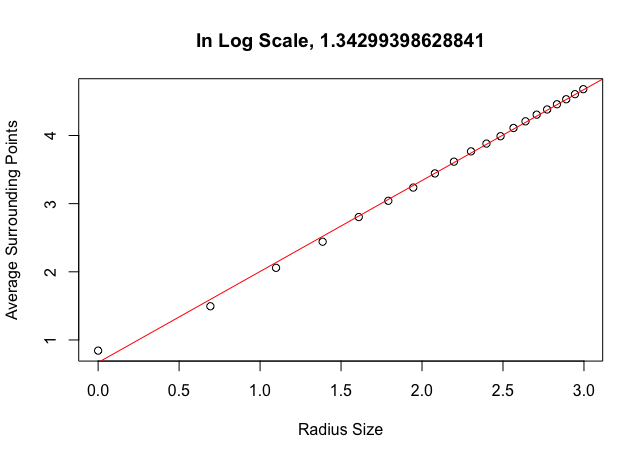}
\end{center}

\subsection{Brownian Earthworm}
We proceeded in a similar fashion with the Brownian earthworm model. However, in this case, we did not know what the actual value for the set of holes as this is presently an open question.

After performing hundreds of simulations with the number of steps ranging from $10$ to $10,000,000$, we performed the Counting Method. We calculated the Hausdorff dimension for the set of holes to be approximately $1.5539$.

\begin{center}
\includegraphics[width=0.6\linewidth]{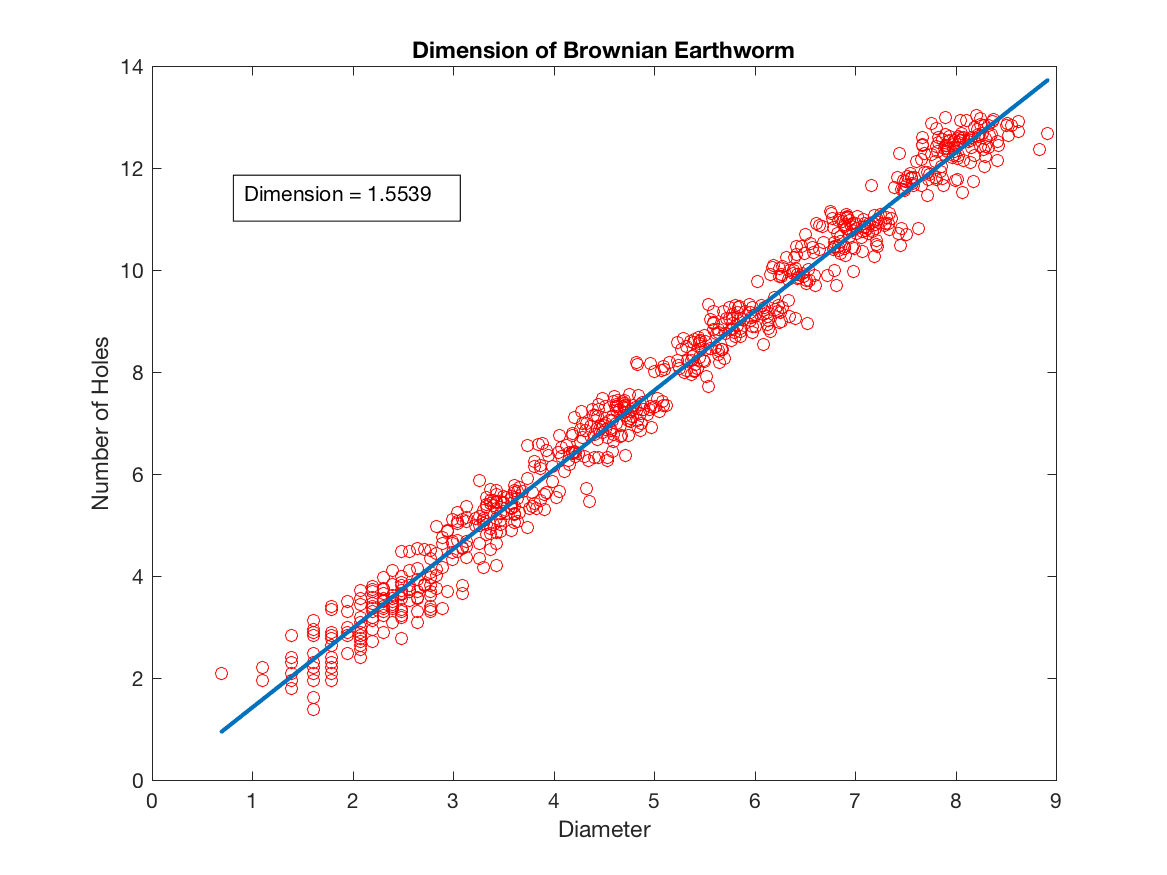}
\end{center}

We again used the Averaging Method to double check our results. This second method gave a fairly similar value for the dimension: $1.5458$. Since our two methods gave similar results, we gained confidence in our methodologies for calculating the Hausdorff dimension.

\begin{center}
\includegraphics[width=0.65\linewidth]{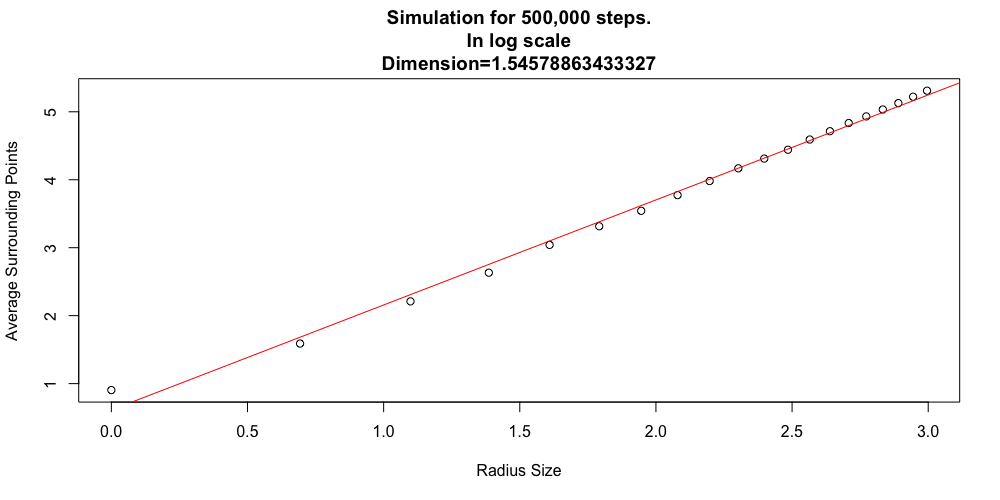}
\end{center}

\section{Conclusions}

After collecting our findings, we consulted with the instructor (Krzysztof Burdzy) about our progress. As it turns out, our methods seem to support the conjecture that the  Hausdorff dimension of the set of holes of the Brownian earthworm is close to, and strictly greater than, $1.5$. The conjecture was based on the simulation results for a closely related model presented in \cite{BHP}.

Toward the end of this project, we began performing some preliminary analysis on the connected components of the Brownian earthworm. In particular, we wanted to consider the growth rate of the number of connected components. Our initial tests seemed to indicate that the plurality of connected components were singletons, although these singletons accounted for only a very small percentage of the total ``area'' of the holes.

%\newpage

%\bibliographystyle{unsrt}
%\bibliography{mybib}

\end{document}